\documentclass{article}
\usepackage{amsfonts}
\usepackage{graphicx}
\usepackage{amsmath}
\usepackage{amssymb}

\begin{document}

\author{Fidel Jos\'{e} Fern\'{a}ndez y Fern\'{a}ndez Arroyo
\and Departamento de Matem\'{a}ticas Fundamentales, UNED
\and C/ Senda del Rey 9, 28040 Madrid, Spain
\and ffernan@mat.uned.es}
\title{The Hahn-Banach Theorem: a proof of the equivalence between the analytic and
geometric versions}
\maketitle

\begin{abstract}
We present here a simple and direct proof of the classic geometric version of
the Hahn-Banach Theorem from its analytic version, in the real case. The reverse
implication, and the direct proofs of both versions, are well known. They are
summarized here for reader's convenience. For the complex case, in both versions the Hahn-Banach Theorem is deduced from the real case, as is well known.
\end{abstract}



\section{Introduction}

The Hahn-Banach Theorem and its applications are essential in Functional
Analysis (see for instance \cite{nb}, \cite{s} and \cite{v}). There are two classical versions, the analytic and the geometric one; both are proved using the Axiom of Choice (see for example \cite{b}). The geometric form also allows to deduce the analytic form, which in its general form (for locally convex spaces) is also called the Separation Theorem. The objective of this work is to present a simple proof of the reverse
implication;\ i.e., to deduce the geometric form from the analytic form, in
the general context of locally convex spaces. The use of a seminorm (which
may not be a norm in general) seems essential. I think that this proof is
of some interest. Several very simple examples are added.

\section{Hahn-Banach Theorem (analytic version, real case)}

Let $E$ be a vector space over $\mathbb{R}$, $p:E\rightarrow\mathbb{R}^{+}$ a seminorm, $L$ a vector subspace of $E$, and $f:L\rightarrow\mathbb{R}$ a linear function verifying $\vert f(x)\vert\leq p(x)$, for every $x\in L$. Then, there exists a linear function $g:E\rightarrow\mathbb{R}$ which extends $f$ (i.e., the restriction of $g$ to $L$ coincides with $f$) and verifies $\vert g(x)\vert\leq p(x)$, for
every $x\in E$.

The direct proof is well known. It uses the Zorn Lemma, which is equivalent to the Axiom of Choice. Nevertheless, we will recall here the scheme of the proof. If we consider the family $A$ of all pairs $(H,h)$, where $H$ is a vector subspace of $E$ which contains $L$  and $h:H\rightarrow\mathbb{R}$ is a linear function which extends $f$ and verifies $\vert h\vert\leq p$ (i.e., $\vert h(x)\vert \leq p(x)$, for every $x\in H$), then this family is not empty ($(L,f)\in{A}$), and we can consider in $A$ the relation $\preceq$ given by $(H_{1},h_{1})\preceq(H_{2},h_{2})$ if $H_{1}\subseteq H_{2}$ and $h_{2}$ extends $h_{1}$. It is easily checked that $\preceq$ is an order relation, and that $({A},\preceq)$ is inductive. The Zorn Lemma guarantees that there exists a maximal element in $({A},\preceq)$, say $(G,g)$. The delicate point is to prove $G=E$, and we shall do this in detail.

If $G\neq E$, then there exists $z\in E\setminus G$. On the other
hand, for every $x_{1},x_{2}\in G$, we get $g(x_{2})-g(x_{1})=g(x_{2}-x_{1})\leq\vert g(x_{2}-x_{1})\vert\leq p(x_{2}-x_{1})=p((x_{2}+z)-(x_{1}+z))\leq p(x_{2}+z)+p(-(x_{1}+z))=p(x_{2}+z)+p(x_{1}+z)$, and therefore $-g(x_{1})-p(x_{1}+z)\leq-g(x_{2})+p(x_{2}+z)$, for every $x_{1},x_{2}\in G$. Let $\gamma=\inf\{-g(x)+p(x+z)\}_{x\in L}$. Then $-g(x_{1})-p(x_{1}+z)\leq\gamma\leq-g(x_{2})+p(x_{2}+z)$, for every $x_{1},x_{2}\in G$. So, for every $x\in G$ and every $t>0$, we get  $-g(\frac{1}{t}x)-p(\frac{1}{t}x+z)\leq\gamma\leq-g(\frac{1}{t}x)+p(\frac{1}{t}x+z)\Rightarrow -g(x)-p(x+tz)\leq t\gamma\leq-g(x)+p(x+tz)$, and $-p(x+tz)\leq g(x)+t\gamma\leq p(x+tz)\Leftrightarrow\vert g(x)+t\gamma\vert \leq p(x+tz)$, for every $x\in G$ and every $t>0$. If $t<0$ and $x\in G$, then $-t>0$, $-x\in G$, and $\vert g(x)+t\gamma\vert=\vert g(x)-(-t)\gamma\vert=\vert-g(-x)-(-t)\gamma\vert=\vert g(-x)+(-t)\gamma\vert\leq p(-x+(-t)z)=p(-(x+tz))=p(x+tz)$. Finally, if $t=0$ and $x\in L$, then $\vert g(x)+t\gamma\vert=\vert g(x)\vert\leq p(x)=p(x+tz)$. So, it is easily checked that, if we define $r:G\oplus\langle z\rangle\rightarrow\mathbb{R}$ by $r(x+tz)=g(x)+t\gamma$, for every $t\in\mathbb{R}$ and every $x\in G$, then $r$ is well defined and it is linear,  $(G\oplus\langle z\rangle,r)\in{A}$, $(G,g)\preceq(G\oplus\langle z\rangle ,r)$, and \ $G\neq G\oplus\langle z\rangle $; which is a contradiction, since $(G,g)$ is maximal in $({A},\preceq)$. We conclude that $G=E$, and it is immediate to see that $(G,g)$ satisfies the required conditions.

Obviously we can prove, in an analogous form, the same version for the case in which the seminorm $p$ is a norm.

\section{Hahn-Banach Theorem (geometric version, re\-al case)}

Let $E$ be a topological vector space over $\mathbb{R}$, $A$ an open and convex subset of $E$, and $S$ a subspace of $E$ such that $S\cap A=\emptyset$. Then, there exists a hyperplane $H$ of $E$ verifying $S\subseteq H$ and $H\cap A=\emptyset$.

The previous statement is also called the Separation Theorem. The direct proof is also well known. It uses again the Zorn Lemma. We will recall here the principal steps. If $A=\emptyset$, then the result is trivial, prolonging a basis of $S$ (let us note that we use also the Zorn Lemma for proving it). So, we will suppose that $A\neq\emptyset$. We consider the family $B$ of all the subspaces $L$ of $E$ such that $S\subseteq L$ and $L\cap A=\emptyset$. $B$ is not empty ($S\in B$), and we can order $B$ by inclusion. It is immediate to see that $(B,\subseteq)$ is inductive. Applying the Zorn Lemma, we get a maximal element $G$ in $(B,\subseteq)$. Let $M=G+\underset{\alpha>0}{\cup}\alpha A$. It is immediate to check that $A\subset M$, $M$ is open, $-M=G+\underset{\alpha<0}{\cup}\alpha A$, and $M\cap(-M)=\emptyset=M\cap G=(-M)\cap G$. We will see that, if $E=G\cup M\cup(-M)$, then $G$ is a hyperplane of $E$ (it is immediate to see that the reverse is also true). Since $A$ is not empty, there exists $a\in A$. If $G$ is not a hyperplane of $E$, then there exists $y\in E\setminus G\oplus\langle a\rangle $; and if $E=G\cup M\cup(-M)$, we can suppose $y\in-M$. The function $f:[0,1]\rightarrow E$ given by $f(t)=ta+(1-t)y$ is continuous, and it is easy to prove that $f^{-1}(G)=\emptyset$  (since $A\cap G=\emptyset$ and $y\notin G\oplus\langle a\rangle $). Moreover, $f(1)=a\in M$, $f(0)=y\in-M$, and so $[0,1]=f^{-1}(M)\cup f^{-1}(-M)$, which is a contradiction since the set $[0,1]$ is connected in $\mathbb{R}$, $f^{-1}(M)$ and $f^{-1}(-M)$ are open disjoint sets, and both are nonempty. Therefore, if $G$ is not a hyperplane of $E$, then $E\neq G\cup M\cup(-M)$. Let $z\in E\setminus G\cup M\cup(-M)$. It is easy to check that $A\cap(G\oplus\langle z\rangle)=\emptyset$, and $G$ is not maximal in $B$. We conclude that $G$ must be a hyperplane of $E$. Moreover, $S\subseteq G$ and $A\cap G=\emptyset$. This completes the classical proof.

It is also well known that we can prove, following the same way, this version for the case in which $E$ is a normed space.

\section{Proof (of the analytic version from the geometric version)}

It is also well known the direct proof of the analytic version from the geometric version.

We will recall here the proof. We suppose the geometric version. Let $E$ be a vector space over $\mathbb{R}$, $p:E\rightarrow\mathbb{R}^{+}$ a seminorm, $L$ a vector subspace of $E$, and $f:L\rightarrow\mathbb{R}$ a linear function verifying $\vert f(x)\vert\leq p(x)$, for every $x\in L$. We want to see that there exists a linear function $g:E\rightarrow\mathbb{R}$ which extends $f$ and verifies $\vert g(x)\vert\leq p(x)$, for every $x\in E$. If $f=0$ (null constant function), then we consider $g=0$, and the result is trivial. If $f\neq0$, then there exists $y\in L$ such that $f(y)=1$. We consider the topology induced by the seminorm $p$. The set $A=\{x\in A \mid p(y-x)<1\}$ is open and convex. On the other hand, $\text{Ker}(f)=\{x\in L\mid f(x)=0\}$ is a vector subspace of $E$, and it is trivial to see that $A\cap\text{Ker}(f)=\emptyset$. From the geometric version of the Hahn-Banach Theorem, there exists a hyperplane $H$ of $E$ verifying $\text{Ker}(f)\subseteq H$ and $H\cap A=\emptyset$. Then, $E=H\oplus\langle y\rangle$. We consider the function $g:E\rightarrow\mathbb{R}$ given by $g(x+ty)=t$ ($x\in H$, $t\in\mathbb{R}$). Obviously, $g$ is well defined and it is linear. It is immediate to check that $g$ extends $f$. We will see that, for every $x\in H$ and every $t\in\mathbb{R}$, $\vert g(x+ty)\vert=\vert t\vert\leq p(x+ty)$. In fact, if $t=0$ then the result is trivial; and if $t\neq0$, then $-\frac{1}{t}x\in H\subseteq E\setminus A$, and therefore $1\leq p(y-(-\frac{1}{t}x))=p(y+\frac{1}{t}x)=p(\frac{1}{t}(ty+x))=\frac{1}{\vert t\vert}p(ty+x)$. So, the linear function $g$ satisfies the required conditions.

It is immediate to see that we can obtain a similar proof for the case of normed spaces ($p$ being a norm).

\section{Results}

We present here a general proof, in the real case, for the reverse
implication. I have not been able to give a similar proof for the case of
normed spaces, and in fact we will see in the second and third examples below that the
seminorm which we obtain cannot be a norm.

\subsection{Proof (of the geometric version from the analytic version)}

Let $E$ be a topological vector space over $\mathbb{R}$, $A$ an open and convex subset of $E$, and $S$ a subspace of $E$ such that $S\cap A=\emptyset$. We want to see that there exists a hyperplane $H$ of $E$ verifying $S\subseteq H$ and $H\cap A=\emptyset$. If $A=\emptyset$, then the result is trivial. So, we will suppose $A\neq\emptyset$. We consider the open set $B=\underset{\alpha>0}{\cup}\alpha A$, which
is obviously nonempty. It is easy to see that $B$ is convex, $B\cap S=\emptyset$, and $\alpha B=B$ for every $\alpha>0$. Let $x\in B$. The set $B-x$ is an open and convex neighborhood of the origin; but it is not balanced. Nevertheless, the set $D=(B-x)\cap(x-B)$ is an open and absolutely balanced (convex and balanced) neighborhood of the origin. Let $p$ be the Minkowski functional of $D$. Then, $p$ is a continuous seminorm, and $p(e)<1$ ($e\in E$) if and only if $e\in D$. We remark that $x\in B\subseteq E\setminus S$. Let $L$ be the linear span of $S\cup\{x\}$. Obviously, for every $y\in L$, there exists a unique $z\in S$, and a unique $t\in\mathbb{R}$, such that $y=z+tx$. We define $f(y)=f(z+tx)=t$. It is immediately checked that the function $f:L\rightarrow\mathbb{R}$ is well defined, it is linear, and $f_{|S}\equiv 0$. For every $z\in S$ and every $t\neq0$, we get $\frac{1}{t}z+x\notin x-B$ (since $B\cap S=\emptyset$). Therefore, for every $z\in S$ and every $t\neq0$, we get $\frac{1}{t}z+x\notin(x-B)\cap(B-x)=D$, and so $p(\frac{1}{t}z+x)\geq 1$. In consequence, $p(z+tx)=p(t(\frac{1}{t}z+x))=\vert t\vert p(\frac{1}{t}z+x)\geq\vert t\vert 1=\vert t\vert =\vert f(z+tx)\vert$. And for every $z\in S$, if $t=0$ we get $\vert f(z)\vert=0\leq p(z)$. It follows that $\vert f(y)\vert $ $\leq p(y)$, for every $y\in L$. From the analytic version of the Hahn-Banach Theorem, in the real case, we deduce that there exists a linear function $g:E\rightarrow\mathbb{R}$ such that $\vert g(e)\vert\leq p(e)$, for every $e\in E$, and the restriction of $g$ to $L$ coincides with $f$. We will see that $H=\text{Ker}(g)$ is a hyperplane containing $S$ and contained in $E\setminus A$. First, since $x\in L$ and $g(x)=f(x)=1$, we get $g\neq0$ , and therefore $H=\text{Ker}(f)$ is a hyperplane. On the other hand, $S\subset L$, and for every $z\in L$ we get $g(z)=f(z)=f(z+0x)=0$; so, $S\subset\text{Ker}(g)=H$. Lastly, we will prove $H\cap A=\emptyset$. In fact, we will see $H\cap B=\emptyset$. If we suppose that there exists $y\in B$ such that $g(y)=0$, then, since the function $g$ is linear, for every $t\in\mathbb{R}$ we get $g(x+ty)=g(x)=1$. Since $\vert g(e)\vert\leq p(e)$, for every $e\in E$,  we get $1=g(x+ty)\leq\vert g(x+ty)\vert\leq p(x+ty)$, and therefore $x+ty\notin D=(B-x)\cap(x-B)$, for every $t\in\mathbb{R}$. Since $y\in B$, if $t<0$ we get $x+ty\in x-B$, for every $t<0$; and since $x+ty\notin D=(B-x)\cap(x-B)$, for every $t\in\mathbb{R}$, it follows that $x+ty\notin B-x$, for every $t<0$. We deduce that $ty\notin B-2x$, for every $t<0$. As $E$ is a topological vector space, the function $h:\mathbb{R}\rightarrow E$ given by $h(y)=ty$ is continuous, and obviously $h(0)=0$. On the other hand, since $x\in B$, we get $2x\in2B=B$, and therefore $0\in B-2x$, where the set $B-2x$ is open. So, $0\in h^{-1}(B-2x)$, being open the set $h^{-1}(B-2x)$; and therefore, there exists $t<0$ such that $h(t)=ty\in B-2x$. This is a contradiction with the previous result $ty\notin B-2x$, for every $t<0$. We conclude that $H\cap B=\emptyset$.

\subsection{Remarks}

\begin{enumerate}
\item As a curiosity, we note that, in the previous proof, $x\notin D$ (since $0\notin B$), and so $p(x)\geq1$. On the other hand, we see that $x\in\overline{D}$. In fact, it is immediate to see that, for every $t\in(-1,1)$, $tx\in(B-x)\cap(x-B)=D$. Since the function $v:\mathbb{R}\rightarrow X$ given by $v(t)=tx$ is continuous, we conclude that $x\in\overline{D}$ and $p(x)\leq1$. So, we get $p(x)=1$.

\item On the other hand, let us note that, if $g:E\rightarrow\mathbb{R}$ is a linear function, and $g$ extends the function $f$ obtained in the previous proof (i.e., $g(z+tx)=t$, \ for every $z\in S$ and every $t\in\mathbb{R}$), then $\vert g\vert\leq p$ if and only if $A\cap\text{Ker}(g)=\emptyset$. We have already proved that if $\vert g\vert\leq p$ (which is equivalent to $g\leq p$, as it is easy to check) then $A\cap\text{Ker}(g)=\emptyset$ (which is equivalent to $B\cap\text{Ker}(g)=\emptyset$). We will see now that, if $A\cap\text{Ker}(g)=\emptyset$, then $\vert g\vert \leq p$. By hypothesis, $S\subseteq\text{Ker}(g)$ and $g(x)=1$; so, $E=\langle x\rangle\oplus\text{Ker}(g)$. For every $e\in E$, there exists $t\in\mathbb{R}$ and $z\in\text{Ker}(g)$ such that $e=tx+z$. Let us note that $g(e)=t$. If $t=0$, then $\vert g(e)\vert=0\leq p(e)$. If $t\neq0$, then $(-\frac{1}{t})z=-\frac{1}{t}z\in\text{Ker}(g)$. Since $A\cap\text{Ker}(g)=\emptyset$, and $\text{Ker}(g)$ is a subspace of $E$, we get $B\cap\text{Ker}(g)=\emptyset$. So, $-\frac{1}{t}z\notin B$. Therefore, $x+\frac{1}{t}z=x-(-\frac{1}{t}z)\notin x-B$, and we get $x+\frac{1}{t}z\notin D$. We deduce that $1\leq p(x+\frac{1}{t}z)=p(\frac{1}{t}(tx+z))=\frac{1}{\vert t\vert}p(tx+z)\Leftrightarrow\vert t\vert=\vert g(e)\vert\leq p(tx+z)=p(e)$. This completes the proof.
\end{enumerate}

\subsection{Examples}

\begin{enumerate}
\item Let $E=\mathbb{R}^{2}$ be the plane with the usual structure of topological vector
space, $A=\{(x,y)\in\mathbb{R}^{2}\mid(x-2)^{2}+y^{2}<2\}$ the open circle with center in $(2,0)$ and radius $\sqrt{2}$, and $S=\{(0,0)\}$. In this case, $B=\underset{\alpha>0}{\cup}\alpha A=\{(x,y)\in\mathbb{R}^{2}\mid x>0,\vert y\vert <x\}$ is the region of the plane on the right of the two lines $y=x$ and $y=-x$, as it is easy to check. In fact, $(x,y)\in B=\underset{\alpha>0}{\cup}\alpha A\Leftrightarrow\exists\alpha>0\mid(x,y)\in\alpha A\Leftrightarrow\exists\alpha>0\mid2\alpha^{2}-4\alpha x+(x^{2}+y^{2})<0$. Since the minimum value of the function $g(\alpha)=2\alpha^{2}-4\alpha x+(x^{2}+y^{2})$ is obtained at $x$, we get $(x,y)\in B\Leftrightarrow x>0$ and $2x^{2}-4x^{2}+x^{2}+y^{2}=-x^{2}+y^{2}<0\Leftrightarrow x>0$ and $\vert y\vert <x$. We choose $\overline{x}\in B$. For instance, let $\overline{x}=(1,0)$. We obtain $B-\overline{x}=\{(x,y)\in\mathbb{R}^{2}\mid x>-1,\vert y\vert <x+1\}$, and $D=(B-\overline{x})\cap(\overline{x}-B)=\{(x,y)\in\mathbb{R}^{2}\mid\vert x\vert<1,\vert y\vert<1-\vert x\vert\}$. Let $p$ be the Minkowski functional associated to $D$. We get, for every $(x,y)\in\mathbb{R}^{2}$, $p(x,y)=\inf\{t>0\mid(x,y)\in tD\}=\inf\{t>0\mid\vert\frac{1}{t}x\vert<1,\vert \frac{1}{t}y\vert <1-\vert\frac{1}{t}x\vert\}=\inf\{t>0\mid\vert x\vert<t,\vert x\vert+\vert y\vert<t\}=\vert x\vert+\vert y\vert$. In this case, $L=S\oplus\langle\overline{x}\rangle=\{(x,y)\in\mathbb{R}^{2}\mid y=0\}$ is the horizontal axis. Moreover, $\forall x\in\mathbb{R}$, $f(x,0)=xf(\overline{x})=x$. Let $g:E\rightarrow\mathbb{R}$ be a linear function verifying $g(x,0)=x$ and $\vert g(x,y)\vert\leq p(x,y)=\vert x\vert+\vert y\vert$, for every $(x,y)\in\mathbb{R}^{2}$ (we know that a such linear function exists by the Hahn-Banach Theorem in its analytic version). Since $g$ is linear, there exist $a,b\in\mathbb{R}$ such that  $g(x,y)=ax+by$, for every $(x,y)\in\mathbb{R}^{2}$. We get, $\forall(x,y)\in\mathbb{R}^{2}$, $g(x,0)=ax+b0=ax=x$ and $\vert g(x,y)\vert=\vert ax+by\vert\leq p(x,y)=\vert x\vert+\vert y\vert\Leftrightarrow a=1$ and $\forall(x,y)\in\mathbb{R}^{2}$, $\vert x+by\vert\leq\vert x\vert+\vert y\vert\Leftrightarrow a=1$ and $\vert b\vert\leq1$. So, $\text{Ker}(g)=\{(x,y)\in\mathbb{R}^{2}\mid x+by=0\}$, with $\vert b\vert\leq1$. It is easy to check that $\text{Ker}(g)$ is a straight line passing through the origin, that it forms an angle between $45^{\circ}$ and $135^{\circ}$ with the horizontal axis, and that $A\cap\text{Ker}(g)=\emptyset$. In fact, such straight lines are all those containing the origin and not meeting $A$.

\item Let $E=\mathbb{R}^{3}$, with the usual topology. Let $S=\{(x,y,z)\in\mathbb{R}^{3}\mid x=y=0\}$, and let $A=\{(x,y,z)\in\mathbb{R}^{3}\mid x>0\}$. Obviously, $S$ is a vector subspace of $E$, $S\cap A=\emptyset$, $A\neq\emptyset$, and the set $A$ is open and convex. In this case, $B=\underset{\alpha>0}{\cup}\alpha A=A$. Let $\overline{x}=(1,-3,0)\in A$. We get $D=(B-\overline{x})\cap(\overline{x}-B)=\{(x,y,z)\in\mathbb{R}^{3}\mid\vert x\vert<1\}$, and the Minkowski functional $\Psi$ associated to $D$ is given by $\Psi(x,y,z)=\vert x\vert$, as it is easily checked. Let us note that, in this example, the seminorm $\Psi$ is not a norm. In this case, $L=S\oplus\langle\overline{x}\rangle=\{(x,y,z)\in\mathbb{R}^3\mid 3x+y=0\}$, and the function $f:L\rightarrow\mathbb{R}$ is given by $f(x,y,z)=x=-\frac{1}{3}y$. We remark that $\vert f(x,y,z)\vert=\Psi(x,y,z)$. Let $g:E\rightarrow\mathbb{R}$ be a linear function verifying $\vert g(x,y,z)\vert\leq\Psi(x,y,z)=\vert x\vert $, for every $(x,y,z)\in E$, and also such that if $x=-\frac{1}{3}y$ then $g(x,y,z)=x$. It is easy to check that $g(x,y,z)=x$, for every $(x,y,z)\in\mathbb{R}^{3}$, and $\text{Ker}(g)=\{(x,y,z)\in\mathbb{R}^{3}\mid x=0\}$ is the only hyperplane containing $S$ and not meeting $A$.

\item Let $E=\mathcal C[0,1]$ be the vector space of continuous real functions defined in $[0,1]$, with the locally convex topology of the punctual convergence. Let $S=\{f\in E\mid f(1)=f(0)=0\}$, and $A=\{f\in E\mid f(1)<0\}$. It is immediately checked that $S$ is a vector subspace of $E$, $S\cap A=\emptyset$, $A\neq\emptyset$, and the set $A$ is open and convex. Moreover, in this case, $B=\underset{\alpha>0}{\cup}\alpha A=A$. Let $h(x)=x-2$, for every $x\in[0,1]$. Obviously, $h\in A$. We get $D=(B-h)\cap(h-B)=\{h\in E\mid\vert h(1)\vert<1\}$, and the Minkowski functional $\Psi$ associated to $D$ is given by $\Psi(h)=\vert h(1)\vert $ ($h\in E$), as it is easily checked. Let us note that this seminorm $\Psi$ is not a norm. In this case, $L=S\oplus\langle h\rangle=\{f\in E\mid f(0)=2f(1)\}$, and the function $F:L\rightarrow\mathbb{R}$ is given by $F(f)=-\frac{1}{2}f(0)=-f(1)$. We remark that $\vert F(f)\vert =\Psi(f)$, for every $f\in G$. Let $G:E\rightarrow\mathbb{R}$ be a linear function verifying $\vert G(f)\vert\leq\Psi(f)=\vert f(1)\vert $, for every $g\in E$, and also such that if $f(0)=2f(1)$ ($f\in E$), then $G(f)=-f(1)$. It is easy to check that $G(f)=-f(1)$, for every $f\in E$, and $\text{Ker}(g)=\{f\in E\mid f(1)=0\}$ is the only hyperplane containing $S$ and not meeting $A$.
\end{enumerate}

\section{Conclusions}

It is possible to find a simple and general proof of the equivalence between
the two classical versions (analytic and geometric) of this very well known
Theorem, in the general context of real topological vector spaces. This proof, which can be used in different examples, shows again the deep connection between both apparently different versions.

The original part of the paper, to the best of the author's knowledge, is the proof of the geometric version from the analytic version. Some examples are shown. In contrast to the other proofs here summarized, I have not found a similar proof of this implication for normed spaces.

\section{Acknowledgements}

I wish to devote this work to Dr. Pedro Jim\'{e}nez Guerra. I am very grateful to Dr. Antonio F\'{e}lix Costa Gonz\'{a}lez, Dr. Jorge L\'{o}pez Abad, Dr. Javier P\'{e}rez \'{A}lvarez, and Dr. Jos\'{e} Carlos Sierra Garc\'{\i}a for their help. I am also very grateful to the referee for his valuable comments.

\end{document}